\documentclass{artjlt}

\usepackage{amsmath,amsfonts,amssymb}
\usepackage{mathrsfs}
\usepackage{bbm}
\usepackage{ stmaryrd }
\usepackage{amsfonts}
\usepackage{multirow}
\usepackage{makecell}

\usepackage{hyperref}
\hypersetup{
	colorlinks=true,
	linkcolor=blue,
	filecolor=magenta,  
	urlcolor=cyan,
}
\usepackage[capitalise]{cleveref}

\def\C{\mathbb{C}}

\def\R{\mathbb{R}}
\def\N{\mathbb{N}}
\def\Z{\mathbb{Z}}

\def\UU{\mathcal{U}}

\def\a{\mathfrak{a}}
\def\b{\mathfrak{b}}
\def\c{\mathfrak{c}}

\def\n{\mathfrak{n}}
\def\p{\mathfrak{p}}
\def\q{\mathfrak{q}}
\def\g{\mathfrak{g}}

\def\t{\mathfrak{t}}
\def\h{\mathfrak{h}}

\def\k{\mathfrak{k}}
\def\l{\mathfrak{l}}
\def\s{\mathfrak{s}}
\def\o{\mathfrak{o}}

\def\ol{\overline}

\def\sub{\subseteq}

\def\Aut{\operatorname{Aut}}

\def\ker{\operatorname{ker}}

\def\dim{\operatorname{dim}}

\def\id{\operatorname{id}}

\def\span{\operatorname{span}}

\def\Ad{\operatorname{Ad}}

\title{Iwasawa Decomposition for Lie Superalgebras}                                     
\author{Alexander Sherman\thanks{The author would like to thank his advisor, Vera Serganova, for suggesting this problem and for helpful discussions along the way.  The author also thanks Shifra Reif for many helpful discussions.  Finally, thank you to an anonymous referee for a thorough reading and many helpful comments and suggestions.  This research was partially supported by NSF grant DMS-1701532.}}                 
\lastname{Iwasawa decomposition}  

\msc{17B22, 17B20, 17B40}    

\keywords{Lie superalgebras, symmetric pairs, root systems}         

\address{%
Alexander Sherman\\               
Dept. of Mathematics\\
University of California at Berkeley\\
Berkeley, CA 94720\\
USA\\            
xandersherm@gmail.com               
}

\begin{document}


\maketitle

\begin{abstract}
Let $\g$ be a basic simple Lie superalgebra over an algebraically closed field of characteristic zero, and $\theta$ an involution of $\g$ preserving a nondegenerate invariant form.  We prove that at least one of $\theta$ or $\delta\circ\theta$ admits an Iwasawa decomposition, where $\delta$ is the canonical grading automorphism $\delta(x)=(-1)^{\ol{x}}x$.  The proof uses the notion of generalized root systems as developed by Serganova, and follows from a more general result on centralizers of certain tori coming from semisimple automorphisms of the Lie superalgebra $\g$.    
\end{abstract}

\section{Introduction}\label{sec_intro}

Let $(\g,\k)$ be a symmetric pair coming from an involution $\theta$ of $\g$, where $\g$ is a reductive Lie algebra over an algebraically closed field of characteristic zero.  Then we have the well-known Iwasawa decomposition of $\g$ given by $\g=\k\oplus\a\oplus\n$, which plays an important role in the study of symmetric spaces.  Here $\a$ is a maximal toral subalgebra of $\p$, where $\p$ is the $(-1)$-eigenspace of $\theta$, and $\n$ is the sum of positive weight spaces for the adjoint action of $\a$ on $\g$, for some choice of positivity.

A close analogue of this situation for Lie superalgebras is to consider a supersymmetric pair $(\g,\k)$ coming from an involution $\theta$ of $\g$, where $\g$ is a basic simple Lie superalgebra, i.e.\ $\g$ is simple, admits a nondegenerate invariant form, and $\g_{\ol{0}}$ is reductive.  However it is well known that even if $\theta$ preserves an invariant form on $\g$, there need not be an Iwasawa decomposition in this setting.  We seek to clarify the situation by proving that if $\theta$ does not admit an Iwasawa decomposition then $\delta\circ\theta$ does, where $\delta(x)=(-1)^{\ol{x}}x$.  Note that $\delta\circ\theta|_{\g_{\ol{0}}}=\theta|_{\g_{\ol{0}}}$, so these involutions are closely related to one another. The theorem fails if we do not assume that $\theta$ preserves a nondegenerate form -- see Remark \ref{counterexample} for an example.

An important consequence of the Iwasawa decomposition is the existence of a Borel subalgebra of $\g$ complementary to $\k$ -- in particular one can find a Borel subalgebra containing $\a\oplus\n$.  Thus a corresponding global symmetric space $G/K$ will be a spherical variety.  Algebraic symmetric spaces give rise to a beautiful and well-understood source of spherical varieties.  Many of the features and structures enjoyed by symmetric varieties have been generalized to spherical varieties such as the little Weyl group (\cite{knop1990weylgruppe} and \cite{knop1994asymptotic}), wonderful compactifications (\cite{de1983complete}), and (restricted) root systems (\cite{brion1990vers} and \cite{knop1996automorphisms}).  The author has begun a study of spherical supervarieties and their properties in \cite{sherman2021spherical} and \cite{sherman2020spherical}, and this paper shows that many symmetric supervarieties are spherical using the existence of an Iwasawa decomposition.  

Another important use of the Iwasawa decomposition is in the study of invariant differential operators on the symmetric space $G/K$.  One uses the decomposition to obtain a natural projection $\UU\g/\UU\g\k\to S(\a)$, giving rise to the Harish-Chandra homomorphism associated to this pair.  In \cite{alldridge2012harish} a characterization of the image of the Harish-Chandra homomorphism was given for supersymmetric pairs which admit an Iwasawa decomposition.  See also \cite{sahi2016capelli} and \cite{sahi2020capelli} for work on the Capelli eigenvalue problem on symmetric supervarieties.  

Restricted root systems coming from supersymmetric pairs were used in \cite{sergeev2004deformed} to construct new families of Calogero-Moser-Sutherland systems which are completely integrable.  In at the end of Section 6, we explain a relationship between the deformed root systems used in \cite{sergeev2004deformed} and the restricted root systems obtained from supersymmetric pairs.

We now explain what will be shown.  Let $V$ be a vector space with a symmetric bilinear form, and $R\sub V\setminus\{0\}$ a finite irreducible generalized reflection root system (GRRS) (see Definition \ref{sec_defs} for full definitions).  GRRSs were defined in \cite{gorelik2017generalized}.  Finite GRRSs are a very mild generalization of generalized root systems (GRSs) as defined in \cite{serganova1996generalizations}, and they are more suitable for our purposes.  An irreducible GRRS should be viewed as the root system of a basic (almost) simple Lie superalgebra $\g$.  

Now let $\theta$ be an automorphism of $V$ preserving both the form and $R$.  This automorphism may come from a semisimple automorphism of $\g$, and if $\theta$ comes from an involution of $\g$ then it will be of order 2.  Write $S\sub R$ for those roots fixed by $\theta$.  A root $\alpha\in R$ is odd if the corresponding root space in $\g$ is odd (for the definition of odd roots in terms of GRRSs, see Definition \ref{defn_real_imag_odd}).  The following theorem is the main technical result upon which all other results are based.
\begin{Theorem}\label{main_grs_thm}
	Let $T\sub S$ be the smallest subset of $S$ containing all odd roots of $S$ and such that if $\alpha\in T,\beta\in S$ with $(\alpha,\beta)\neq0$, then $\beta\in T$.  Then we have one of the following possibilities for $T$:
	\begin{itemize}
		\item $T=\emptyset$;
		\item $T=\{\pm\alpha\}$ for an isotropic root $\alpha$;
		\item $T\sub\span(T)$ is a finite irreducible GRRS containing at least one odd root.
	\end{itemize}
\end{Theorem}

Now either let $\g$ be a basic simple Lie superalgebra not equal to $\p\s\l(2|2)$ or let $\g$ be $\g\l(m|n)$.  Recall that being basic means there is an even invariant nondegenerate form on $\g$.  Let $\theta\in\Aut(\g)$ be a semisimple automorphism preserving such a form.  Let $\h$ be a $\theta$-stable Cartan subalgebra of $\g_{\ol{0}}$.  Then $\theta$ induces an automorphism of the GRRS $R\sub\h^*$ corresponding to the choice of $\h$.  Write $\a\sub\h$ for the sum of the eigenspaces of $\theta$ on $\h$ with eigenvalue not equal to one.  If we write $S$ for the roots fixed by $\theta$, then the centralizer of $\a$ is given by $\c(\a)=\h+\bigoplus\limits_{\alpha\in S}\g_{\alpha}$.  Using Theorem \ref{main_grs_thm} we obtain:

\begin{Theorem}\label{main_thm}
	The Lie superalgebra $\c(\a)$ is an extension of an abelian Lie superalgebra by the product of ideals $\a\times\tilde{\l}\times\l$, where $\l$ is an even semisimple Lie algebra and $\tilde{\l}$ is isomorphic to either a basic simple Lie superalgebra, $\s\l(n|n)$ for some $n\geq 1$, or is trivial.  
\end{Theorem}

Note that if $\g$ is Kac-Moody, then $\c(\a)$ will also be Kac-Moody, see Theorem \ref{centralizater}.  We emphasize that the nontrivial statement in Theorem \ref{main_thm} is that the centralizer has only one simple superalgebra appearing which is not purely even.  This need not be true for centralizers of an arbitrary torus in $\g$ -- in particular it is false for many Levi subalgebras.

In the case when $\theta$ is of order $2$, we can construct $\h$ so that $\a$ is a maximal toral subspace of $\p$, the $(-1)$-eigenspace of $\theta$ acting on $\g$.  Classically it is known that $\c(\a)_{\ol{0}}\sub\a+\k$.  However it is possible that $\c(\a)_{\ol{1}}\cap\p\neq 0$, in which case the Iwasawa decomposition doesn't hold. However Theorem \ref{main_thm} implies that if $\c(\a)_{\ol{1}}\cap\p\neq 0$, then $\c(\a)_{\ol{1}}\sub\p$.  Therefore if we look at $\delta\circ\theta$ instead, where $\delta(x)=(-1)^{\ol{x}}x$ is the canonical grading automorphism, then for this automorphism we have $\c(\a)_{\ol{1}}\sub\k$, and thus the Iwasawa decomposition will hold.  We state this as the following result (where the case of $\p\s\l(2|2)$ is checked separately).

\begin{Theorem}\label{main_thm_iwasawa}
	If $\theta$ is an involution on a basic simple superalgebra or $\g\l(m|n)$ which preserves the nondegenerate invariant form, then either $\theta$ or $\delta\circ\theta$ admits an Iwasawa decomposition.  In particular, either the fixed points of $\theta$ or the fixed points of $\delta\circ\theta$ have a complementary Borel subalgebra.	
\end{Theorem} 

Observe that it is possible for both $\theta$ and $\delta\circ\theta$ to admit Iwasawa decompositions; indeed, in many cases these involutions are conjugate to one another, for example any involution of $\a\b(1|3)$ satisfies this.

Finally, let us summarize the contents of the paper.  In Section 2 we recall the definition of finite GRRSs, state the classification of finite irreducible GRRSs, and prove a few facts we will need later on about them.  In Section 3 we introduce automorphisms of GRRSs and prove Theorem \ref{main_grs_thm}.  Section 4 applies the results from Section 3 to prove Theorem \ref{main_thm}.  Section 5 proves Theorem \ref{main_thm_iwasawa} and describes supersymmetric pairs that fit into our framework.  In Section 6 we study restricted root systems that arise from supersymmetric pairs, discuss their general properties, and relate them to the work of Sergeev and Veselov in \cite{sergeev2004deformed}.  Finally, the appendix classifies all supersymmetric pairs of $\a\g(1|2)$ and $\a\b(1|3)$.

\section{Generalized Reflection Root Systems}\label{sec_defs}

We work over an algebraically closed field $\Bbbk$ of characteristic zero. In \cite{serganova1996generalizations} the notion of a generalized root system (GRS) was introduced, and GRSs were completely classified.  In \cite{gorelik2017generalized}, this notion was generalized to that of a generalized reflection root system (GRRS) that was designed to encompass root systems of affine Lie superalgebras.  Finite GRRSs come from root systems of certain (almost) simple Lie superalgebras, and we have found they are a natural object to consider for our problem.  

The proofs of properties of GRSs stated in \cite{serganova1996generalizations} carry over almost entirely to finite GRRSs.  We will restate some of these results without proof with this understanding.

\begin{Definition}\label{GRRS_def}
	Let $V$ be a finite-dimensional $\Bbbk$-vector space equipped with a symmetric bilinear form $(\cdot,\cdot)$ (not necessarily nondegenerate).  A finite generalized reflection root system (GRRS) is a nonempty finite set $R\sub V\setminus\{0\}$ satisfying the following axioms:
	\begin{enumerate}
		\item  $\span(R)=V$;
		
		\item for $\alpha\in R$, $(\alpha,-)\neq 0$ as an element of $V^*$.
		
		\item for $\alpha,\beta\in R$ with $(\alpha,\alpha)\neq 0$ we have $k_{\alpha,\beta}:=\frac{2(\alpha,\beta)}{(\alpha,\alpha)}\in\Z$ and $r_{\alpha}(\beta):=\beta-k_{\alpha,\beta}\alpha\in R$;
		
		\item for $\alpha\in R$ such that $(\alpha,\alpha)=0$ there exists a bijection $r_{\alpha}:R\to R$ such that $r_{\alpha}(\beta)=\beta$ if $(\alpha,\beta)=0$, and $r_{\alpha}(\beta)=\beta\pm\alpha$ if $(\alpha,\beta)\neq0$;
		
		\item $R=-R$.
	\end{enumerate}
	We call the elements of $R$ roots.
\end{Definition}
For the rest of this paper we will call a finite GRRS $R$ just a GRRS with the understanding that it is finite.  We will not consider infinite GRRSs.
\begin{remark}
	\begin{itemize}
		\item A GRS, as defined in \cite{serganova1996generalizations}, is exactly a GRRS in which the form $(-,-)$ is assumed to be nondegenerate.
		\item We note that (2) is equivalent to saying that for all $\alpha\in R$ the bijection $r_{\alpha}:R\to R$ is nontrivial.  
		\item Another notion of a GRS was given in definition 7.1 in \cite{serganova1996generalizations}.  If one defines $\alpha^\vee=\frac{2}{(\alpha,\alpha)}(\alpha,-)$ for a non-isotropic root $\alpha$ and $\alpha^\vee=(\alpha,-)$ for an isotropic root $\alpha$, then a GRRS is a GRS in the sense of definition 7.1 of \cite{serganova1996generalizations} if and only if $\alpha^\vee\neq\beta^\vee$ for all odd isotropic roots $\alpha,\beta$.  We will see this is the case for all irreducible GRRSs except for $\tilde{A}(1,1)$, which is defined below.
	\end{itemize}
\end{remark}

\begin{Lemma}\label{sub_grs}
	Let $R\sub V$ be a GRRS and suppose $S\sub R$ is a subset of $R$ such that 
	\begin{itemize}
		\item $S=-S$;
		\item for each $\alpha\in S$ there exists $\beta\in S$ such that $(\alpha,\beta)\neq0$;
		\item for each $\alpha\in S$, $r_{\alpha}(S)=S$.
	\end{itemize}
	Then $S\sub\span(S)$ is a GRRS.
\end{Lemma}
\begin{proof}
	This follows from the definition.
\end{proof}

\begin{Definition}\label{defn_real_imag_odd}
	If $R$ is a GRRS we define the subset of real (non-isotropic) and imaginary (isotropic) roots as
	\[
	R_{re}=\{\alpha\in R:(\alpha,\alpha)\neq0\}  \ \ \ \ \ R_{im}=\{\alpha\in R:(\alpha,\alpha)=0\}.
	\]
	Further, we call $\alpha\in R$ odd if $\alpha\in R_{im}$ or $2\alpha\in R_{re}$.  Otherwise we say a root is even.
\end{Definition}

By Chapter VI of \cite{bourbakilie}, $R_{re}\sub\span(R_{re})=U$ will be a (potentially non-reduced) root system in the usual sense, and in particular the form is nondegenerate when restricted to $U$.  Thus we can decompose $U$ as $U=V_1\oplus\cdots\oplus V_k$, where $R_{re}^i:=R_{re}\cap V_i\sub V_i$ is irreducible and $R_{re}=\coprod\limits_i R_{re}^i$.  Let $W_i$ denote the Weyl group of $R_{re}^i$, and let $W=W_1\times\cdots\times W_k$, the Weyl group of $R_{re}\sub U$.  Then $W$ acts naturally on $V$ and preserves $R$ and the form $(-,-)$.  Finally let $V_0$ be the orthogonal complement to  $U$ in $V$ so that
\[
V=V_0\oplus V_1\oplus\dots\oplus V_k,
\]
where $R_{re}\cap V_0=\emptyset$.  We write $p_i:V\to V_i$ $i=0,1,\dots,k$ for the projection maps.   Note that $(-,-)$ may be degenerate when restricted to $V_0$.

A GRRS $R$ is reducible if we can write $R=R'\coprod R''$, where $R'$ and $R''$ are nonempty and orthogonal to one another.  In this case each of $R'$ and $R''$ will form GRRSs in the respective subspaces they span.  A GRRS $R$ is irreducible if it is not reducible.  Every GRRS can be decomposed into a finite direct sum of irreducible GRRSs.

\begin{Proposition}[Prop. 2.6, \cite{serganova1996generalizations}]\label{components_prop}
	For an irreducible GRRS $R$, either $\dim V_0=1$ and $k\leq 2$, or $\dim V_0=0$ and $k\leq 3$.  If $V_0\neq0$, then $p_0(R_{im})=\{\pm v\}$ for some nonzero vector $v\in V_0$.
\end{Proposition}

\begin{remark}
	Proposition \ref{components_prop} in particular implies that if $V_0=0$ then $(-,-)$ is nondegenerate.  If $V_0\neq0$ then $(-,-)$ is degenerate if and only if it restricts to the zero form on $V_0$.
\end{remark}

For the irreducible root system  $R^i_{re}\sub V_i$, we write 
\[
P_i=\{x\in V_i:\frac{2(x,\alpha)}{(\alpha,\alpha)}\in\Z\text{ for all }\alpha\in R_{re}^i\}
\]
for the weight lattice of $V_i$.  

\begin{Definition}
	A $W_i$-orbit $X\sub P_i$ is small if $x-y\in R_{re}^i$ for any $x,y\in X$, where $x\neq\pm y$.  
\end{Definition}

\begin{Proposition}[Prop. 3.5 of \cite{serganova1996generalizations}]\label{projection_non_isotropic}    
	Let $R$ be a GRRS.  Then $p_i(R_{im})$ is a subset of $P_i\setminus\{0\}$, and is the union of small $W_i$-orbits.  In particular $(p_{i}(\alpha),p_i(\alpha))\neq0$ for all $\alpha\in R_{im}$ and $i>0$.
\end{Proposition}
\begin{remark}
	Note that the second statement of Lemma \ref{projection_non_isotropic} follows from Cor. 1.7 of \cite{serganova1996generalizations}.	
\end{remark}

Let $R$ be a GRRS.  Then $R_{im}$ is $W$-invariant, and thus we may break it up into its orbits
\[
R_{im}=R_{im}^1\sqcup\cdots\sqcup R_{im}^m.
\]
We call the orbits the imaginary components of $R$. 

\begin{Lemma}\label{lemma_1}
	Let $R$ be an irreducible GRRS.  If $\alpha,\beta$ are isotropic roots that lie in the same imaginary component of $R$, and $p_i(\alpha)=\pm p_i(\beta)$ for all $i$, then either $\alpha=\pm\beta$ or $\alpha\pm\beta=2p_\ell(\alpha)\in R_{re}^\ell$ for some $\ell\in\{1,\dots,k\}$.
\end{Lemma}

\begin{proof}
	For ease of notation, for a vector $v\in V$ write $v^2:=(v,v)$, and write $p_i(\beta)=\epsilon_ip_i(\alpha)$, where $\epsilon_i=\pm1$.  Then by assumption we have that 
	\[
	0=(\alpha,\alpha)=\sum\limits_ip_i(\alpha)^2.
	\]

	Suppose that $\alpha\neq\pm\beta$.  Since there are at most three terms in the above sum, there must be an $\ell$ such that $\epsilon_\ell$ is distinct from $\epsilon_i$ for all $i\neq \ell$.   We see that in this notation, 
	\[
	(\alpha,\beta)=\sum\limits_i\epsilon_ip_i(\alpha)^2.
	\]  
	If this quantity is 0, then we may add it to $\epsilon_\ell(\alpha,\alpha)$ and find that $2\epsilon_\ell p_\ell(\alpha)^2=0$, hence $p_\ell(\alpha)^2=0$.  However this contradicts Lemma \ref{projection_non_isotropic}.  So we must instead have $(\alpha,\beta)\neq0$, so that by axiom (2) of a GRS, either $\alpha+\beta$ or $\alpha-\beta$ is a root.  It must be real in either case, and therefore cannot have a component in $V_0$ and can only have a nonzero component in one $V_i$ for some $i>0$.  It now follows whichever of $\alpha\pm\beta$ is a root, it will be equal to $2p_i(\alpha)$ for some $i>0$, and we are done.
\end{proof} 

Thm. 5.10 of \cite{serganova1996generalizations} classified irreducible GRSs.  However from an analysis of the proof one sees that it also classifies GRRSs, and only one extra family of GRRSs arises that are not already GRSs, and this is the family $\tilde{A}(n,n)$.  This is verified in \cite{gorelik2017generalized} as well.  In terms of Lie superalgebras, $\tilde{A}(n-1,n-1)$ is the root system of $\p\g\l(n|n)=\g\l(n|n)/\Bbbk I_{n|n}$.  To be precise, if we write $\h\sub\g\l(n|n)$ for the subalgebra of diagonal matrices, then $\h^*$ has a nondegenerate inner product from the supertrace form.  If we take the subspace of $\h^*$ spanned by roots of $\g\l(n|n)$ and restrict the form to it, we get the GRRS $\tilde{A}(n-1,n-1)$.

In the following theorem we give the classification of irreducible GRRSs.  In each case we will describe $R_{re}$ and $R_{im}$.  We will write $W$ for the Weyl group of $R_{re}$ in each case, and $\omega_{i}^{(j)}\in V_j$ for the $i$th fundamental weight of $R_{re}^j$; for instance if $R_{re}^{(2)}=A_{n}$, then $\omega_1^{(2)}\in V_2$ denotes the first fundamental weight of the root system $A_{n}$, i.e.\ the dominant weight corresponding to the standard representation.  In the case that $V_0\neq0$, we write $v\in V_0$ for the element describe in Proposition \ref{components_prop}.

\begin{Theorem}
	The irreducible GRRSs with $R_{im}\neq0$ are as follows.  
	\begin{enumerate}
		\item[(0)] $\tilde{A}(n,n)$, $n\geq 1$: $R_{re}=A_n\sqcup A_n$, $R_{im}=(W\omega_{1}+v)\sqcup(W\omega_n-v)$;
		
		\item $A(0,n)$, $n\geq 1$: $R_{re}=A_n$, $R_{im}=(W\omega_1+v)\sqcup(W\omega_n-v)$;
		
		\item $C(0,n)$, $n\geq 2$: $R_{re}=C_n$, $R_{im}=(W\omega_1+v)\sqcup(W\omega_1-v)$;
		
		\item $A(m,n)$, $m\neq n, m\geq 1$: $R_{re}^1=A_m$, $R_{re}^2=A_n$, $R_{im}=(W(\omega_1^{(1)}+\omega_n^{(2)})+v)\sqcup(W(\omega_m^{(1)}+\omega_{1}^{(2)})-v)$;
		
		\item $A(n,n)$, $n\geq 2$: $R_{re}^1=A_n$, $R_{re}^2=A_n$, $R_{im}=W(\omega_1^{(1)}+\omega_{n}^{(2)})\sqcup W(\omega_n^{(1)}+\omega_1^{(2)})$;
		
		\item $B(m,n)$, $m,n\geq 1$: $R_{re}^1=B_m$, $R_{re}^2=BC_n$, $R_{im}=W(\omega_1^{(1)}+\omega_1^{(2)})$;
		
		\item $G(1,2)$: $R_{re}^1=BC_1$, $R_{re}^2=G_2$, $R_{im}=W(\omega_1^{(1)}+\omega_1^{(2)})$;
		
		\item $D(m,n)$, $m>2$, $n\geq1$: $R_{re}^1=D_m$, $R_{re}^2=C_n$, $R_{im}=W(\omega_1^{(1)}+\omega_1^{(2)})$;
		
		\item $AB(1,3)$: $R_{re}^1=A_1$, $R_{re}^2=B_3$, $R_{im}=W(\omega_1^{(1)}+\omega_3^{(2)})$;
		
		\item $D(2,n)$, $n\geq 1$: $R_{re}^1=A_1$, $R_{re}^2=A_1$, $R_{re}^3=C_n$, $R_{im}=W(\omega_1^{(1)}+\omega_{1}^{(2)}+\omega_{1}^{(3)})$;
		
		\item $D(2,1;a)$: $R_{re}^1=A_1$, $R_{re}^2=A_1$, $R_{re}^3=A_1$, $R_{im}=W(\omega_1^{(1)}+\omega_{1}^{(2)}+\omega_{1}^{(3)})$.
	\end{enumerate}  The only GRRS which is not a GRS (i.e.\ for which the inner product is degenerate) is $\tilde{A}(n,n)$.   
	
	Each inner product is determined up to proportionality, except for $D(2,1;a)$ where we get a family of distinct inner products parametrized by $a\in\C\setminus\{0,-1\}$ modulo an action of $S_3$. Further the inner products on two distinct real components of $D(2,1;a)$ agree if and only if $D(2,1;a)\cong D(2,1)$, which is when $a=1,-2$, or $-1/2$.
\end{Theorem}

\begin{remark}
	The cases (1)-(10) are each the root system of a unique basic simple Lie superalgebra.  The only basic simple Lie superalgebra that is left out in the above classification is $\p\s\l(2|2)$.  This is due to having root spaces of dimension bigger than one.  However using GRRSs we do get $\tilde{A}(1,1)$, which as already stated corresponds to $\p\g\l(2|2)$, whose derived subalgebra is $\p\s\l(2|2)$.
\end{remark}

\begin{Corollary}\label{lin_combo_cor}
	let $\alpha,\beta$ be linearly independent isotropic roots in an irreducible GRRS $R$.  Then for some $i>0$, one of two things must occur:
	\begin{enumerate}
		\item $p_i(\alpha)$ and $p_i(\beta)$ are orthogonal and either $p_i(\alpha)+p_i(\beta)\in R_{re}^i$ or $p_i(\alpha)-p_i(\beta)\in R_{re}^i$;
		\item $2p_i(\alpha)=\pm 2p_i(\beta)\in R_{re}^i$.
	\end{enumerate}
\end{Corollary}  

\begin{proof}
	If $\alpha$ and $\beta$ lie in the same imaginary component of $R$, then $p_i(\alpha)$ and $p_i(\beta)$ lie in the same small $W_i$-orbit. Let $i$ be such that $R_{re}^i$ is one of $A_n,B_n,C_n,$ or $D_n$ and $p_i(R_{im})\sub\pm W\omega_{1}^{(i)}$.  Observe that for these root systems, if $\lambda,\mu\in W\omega_{1}$ then either $\lambda=\pm\mu$ or $\lambda$ is orthogonal to $\mu$.
	
	Now if $p_i(\alpha)\neq\pm p_i(\beta)$ for some $i$, then $p_i(\alpha)$ is orthogonal to $p_i(\beta)$ and by Lemma \ref{projection_non_isotropic} $p_i(\alpha)-p_i(\beta)\in R_{re}^i$ so we are done.  Otherwise, we are in the situation of Lemma \ref{lemma_1}, giving $2p_i(\alpha)=\pm2p_i(\beta)\in R_{re}^i$ for some $i$, and we are done.
	
	If $\alpha$ and $\beta$ lie in distinct imaginary components, then we have $R$ is one of the GRRSs listed in (0)-(4) above.	 But we see that in each case there are two imaginary components and they are swapped under negation.  Thus $\alpha$ and $-\beta$ are in the same imaginary component, so we may apply the argument just given to finish the proof.
\end{proof} 

\section{Automorphisms of weak generalized root systems}

Let $R\sub V$ be an irreducible GRRS and $\theta$ an automorphism of $R$, meaning that $\theta:V\to V$ is a linear isomorphism preserving the bilinear form, with $\theta(R)=R$. Write $S\sub R$ for the roots fixed by $\theta$.  By linearity, we have that $S=-S$, and if $\alpha,\beta\in S$ with $\alpha+\beta\in R$, then $\alpha+\beta\in S$.  We now prove the main technical result of the paper.

\begin{Proposition}\label{odd_connectedness}
	Let $\alpha,\beta$ be linearly independent odd roots of $S$.  Then there exists a real root $\gamma\in R_{re}$ with $\theta(\gamma)=\gamma$ (i.e.\ $\gamma\in S$) such that $(\gamma,\alpha)\neq0$ and $(\gamma,\beta)\neq0$.
\end{Proposition}

\begin{proof}
	We break the proof up into two cases.
	
	\
	
	\emph{Case 1: $\alpha,\beta$ are isotropic}: 
	
	In general, $\theta$ will either preserve all real components $R_{re}^i$ or will permute them in a nontrivial way.  We first deal with the latter case.  If $\theta$ permutes $R_{re}^i$ and $R_{re}^j$, then in particular these root systems must be isomorphic.  Looking at our list, this leaves only (0), (4), (9), and (10) as possibilities.  However, in the cases of (0) and (4) the inner product on each factor of $A_n$ is negative the other, so no such $\theta$ can exist that permutes them.  Further, in the case of (10) such a permutation could only exist if two of the underlying real root systems are isomorphic, i.e.\ their inner products agree, which would give $D(2,1)$.  So it remains to deal with case (9).  
	
	For the case of (9), we may assume that $R_{re}^3$ is preserved by $\theta$.  If $p_3\alpha\neq \pm p_3\beta$ then necessarily $p_3\alpha$ and $p_3\beta$ are orthogonal because they lie in the orbit of $\omega^{(3)}_1$.  By smallness of the orbit of $\omega_1$ in $C_n$ we will have $\gamma=p_3\alpha-p_3\beta\in R_{re}^3$ is fixed by $\theta$, and this will not be orthogonal to $\alpha$ or $\beta$ so that $(\gamma,\alpha)\neq0$ and $(\gamma,\beta)\neq0$. If $p_3\alpha=\pm p_3\beta$ then $\gamma=2p_3\alpha\in R_{re}^3$ works.
	
	If instead $\theta$ preserves each $R_{re}^i$, then each $p_i\alpha$ is fixed by $\theta$ since $\theta\alpha=\alpha$.  We then apply Corollary \ref{lin_combo_cor} to get that there exists an $i$ such that some linear combination of $p_i(\alpha)$ and $p_i(\beta)$ is in $R_{re}^i$ which is not orthogonal to $\alpha$ or $\beta$ and is fixed by $\theta$. 	
	
	\
	
	\emph{Case 2: one of $\alpha,\beta$ non-isotropic}
	
	If $\alpha$ is non-isotropic, then one real component of $R$ must be $BC_n$ for some $n$, hence either $R=G(1,2)$ or $R=B(m,n)$.  If $R=G(1,2)$, then $\alpha=\pm\omega_1^{(1)}$.  Hence if $\beta$ is isotropic then $(p_1(\beta),\alpha)\neq0$ so we can take $\gamma=\alpha$.  If $\beta$ is non-isotropic then $\beta=\pm\omega_1^{(1)}$ as well, so clearly $(\alpha,\beta)\neq0$ and we can again take $\gamma=\alpha$.
	
	If $R=B(m,n)$ and $\beta$ is isotropic, then $p_2\beta=\sigma\omega_1^{(2)}$ for some $\sigma$ in the Weyl group of $BC_n$.  Hence either $p_2\beta=\pm\alpha$, in which case we can take $\gamma=\alpha$, otherwise $\gamma=p_2\beta+\alpha\in BC_n$ works.  If $\beta$ is non-isotropic then either $\beta=\pm\alpha$, in which case we take $\gamma=\alpha$, and otherwise $\gamma=\beta+\alpha\in BC_n$ works.  
\end{proof}

\begin{Corollary}\label{cor_grrs}
	If $S$ either contains at least 2 linearly independent odd roots or no odd roots at all, then $S\sub\span(S)$ is a GRRS. 
\end{Corollary}

\begin{proof}
	We may apply Section \ref{sub_grs} along with Proposition \ref{odd_connectedness} to obtain the result.  
\end{proof}
\begin{remark}
	Note that we could have $S=\{\pm\alpha\}$ for an isotropic root $\alpha$.  For example if we consider $A(0,2)$, the automorphism given by a simple reflection of the Weyl group of $A_2$ will give rise to such a situation.
\end{remark}

Now let $T\sub S$ be the smallest subset of $S$ satisfying:
\begin{enumerate}
	\item all odd roots of $S$ lie in $T$;
	\item if $\alpha\in T$, $\beta\in S$ with $(\alpha,\beta)\neq 0$, then $\beta\in T$. 
\end{enumerate}   
Then $T$ will be orthogonal to $T':=S\setminus T$, and $T'$ will consist of only even roots.

\begin{Proposition}\label{prop_about_T}
	$T'\sub\span(T')$ is a reduced root system.  Further, we have the following possibilities for $T$:
	\begin{enumerate}
		\item $T=\emptyset$.
		\item $T=\{\pm\alpha\}$ for an isotropic root $\alpha$.
		\item $T\sub\span(T)$ is an irreducible GRRS containing at least one odd root.
	\end{enumerate}
	In all cases, $T$ is orthogonal to $T'$ and we have both $S\cap\span(T)=T$ and $S\cap\span(T')=T'$.
\end{Proposition}
\begin{proof}
	The first statement is clear.  For the second statement, if $S\cap R_{im}=\{\pm\alpha\}$ for some $\alpha$, then we claim $T=\{\pm\alpha\}$.  This is because if not then there exists $\beta\in T\setminus\{\pm\alpha\}$ such that $\beta$ must is real and $(\alpha,\beta)\neq 0$.  Thus $r_\beta\alpha$ would be another isotropic root in $T$.  
	
	If $S\cap R_{im}\neq\{\pm\alpha\}$ for some $\alpha$ then either it is empty, or contains two linearly independent isotropic roots.  In the former case $T$ will either be empty or a non-reduced root system which is irreducible (by Proposition \ref{odd_connectedness}) and thus is $BC_n$.  In the latter case $T\sub\span(T)$ is an irreducible GRRS with $T_{im}\neq\emptyset$ by Proposition \ref{odd_connectedness} and Section \ref{sub_grs}.
	
	Now for each possibility of $T$ we always have that the span of the odd roots is equal to the span of all of $T$, as this is true for any irreducible GRRS admitting at least one odd root.  It follows that $\span(T')$ is orthogonal to $\span(T)$.  Since the inner product restricted to $\span(T')$ will be nondegenerate we must have $S\cap\span(T')=T'$.  On the other hand if $\alpha\in T'\cap\span(T)$ we would have that $\alpha$ is an even null vector, a contradiction.
\end{proof}


\begin{Corollary}\label{S_almost_GRRS}
	Either $S\sub\span(S)$ is a GRRS or $S=T'\sqcup\{\pm\alpha\}$ where $T'\sub\span(T')$ is an even reduced root system and $\alpha$ is an isotropic root orthogonal to $T'$.
\end{Corollary}
\section{Application to centralizers of certain tori}

\begin{Lemma}\label{auto_semisimple}
	Suppose that $\g$ is a Lie superalgebra such that: 
	\begin{enumerate}
		\item $\g_{\ol{0}}$ is reductive and $\g_{\ol{1}}$ is a semisimple $\g_{\ol{0}}$-module;
		\item If $\h  \sub\g_{\ol{0}}$ is a Cartan subalgebra (CSA) of $\g_{\ol{0}}$, then it is self-centralizing in $\g$.
		\item For any root $\alpha$ we have $\dim\g_{\alpha}\leq 1$.
	\end{enumerate} Then $\theta\in\Aut(\g)$ is semisimple if and only if $\theta|_{\g_{\ol{0}}}$ is semisimple.  In particular, $\theta$ is semisimple if and only if it preserves a Cartan subalgebra of $\g_{\ol{0}}$.
\end{Lemma}  
\begin{remark}
	Property (2) is equivalent to asking that for any root decomposition of $\g$, each weight space (including the trivial weight space) is of pure parity.  
\end{remark}
\begin{proof}
	By \cite{borel1955semi}, an automorphism of a reductive Lie algebra is semisimple if and only if it preserves a Cartan subalgebra.  Therefore if $\theta|_{\g_{\ol{0}}}$ is semisimple, it preserves a Cartan subalgebra $\h\sub\g_{\ol{0}}$, and thus must act by a permutation on the roots.  Since the root spaces are one-dimensional, it follows that some power of $\theta$ must act by a scalar on each weight space, and thus $\theta$ must be semisimple.
\end{proof} 

Suppose that $\g$ either is a basic simple Lie superalgebra not equal to $\p\s\l(2|2)$ or is $\g\l(m|n)$ for some $m,n$ so that $\g$ satisfies the hypothesis of Lemma \ref{auto_semisimple}.  Let $\theta\in\Aut(\g)$ be a semisimple automorphism of $\g$ which preserves a nondegenerate invariant form on $\g$.  We get an orthogonal decomposition $\g=\k\oplus\p$ of super vector spaces, where $\k$ is the fixed subalgebra of $\theta$, and $\p$ is the sum of the eigenspaces of $\theta$ with nonzero eigenvalues. 

\begin{remark}\label{killing_form_pres}
	The Killing form is nondegenerate for $\s\l(m|n)$ with $m\neq n$, $\o\s\p(m|2n)$ when $m-2n\neq2$ and $m+2n\geq 2$, and for $G(1,2)$ and $AB(1,3)$.  Thus every automorphism of these superalgebras necessarily preserves the form.
\end{remark}

Now suppose $\h\sub\g_{\ol{0}}$ is a Cartan subalgebra which is $\theta$-invariant.  Write $\h=\t\oplus\a$, where $\t=\k\cap\h$ and $\a=\p\cap\h$.  Then $\theta$ induces an automorphism of $\h^*$ preserving the set of roots, $R$, and thus induces an automorphism of the GRRS $R\sub V=\span(R)$.  In the case of $\g\l(m|n)$, $R\sub\span(R)$ will either be $A(m-1,n-1)$ if $m\neq n$ or $\tilde{A}(n-1,n-1)$ if $m=n\neq 1$, and this is the GRRS we consider.  If $m=n=1$, we do not obtain a GRRS but the following will be easy to check in this case anyway.

We keep the notations as above for $S$, $T,$ and $T'$.   Write $\c(\a)$ for the centralizer of $\a$ in $\g$.  Notice that we have $\c(\a)=\h+\bigoplus\limits_{\alpha\in S}\g_\alpha$

\begin{Proposition}\label{centralizer_prop}
	Define the following subalgebras of $\g$:
	\begin{itemize}
		\item $\l$ the subalgebra of $\g$ generated by $\{e_{\alpha}:\alpha\in T'\}$;
		\item $\tilde{\l}$ the subalgebra of $\g$ generated by $\{e_{\alpha}:\alpha\in T\}$.
	\end{itemize} 
	Then $\l$ is a semisimple Lie algebra, and $\tilde{\l}$ either is isomorphic to a basic simple Lie superalgebra, isomorphic to $\s\l(n|n)$ for some $n\geq 1$, or is trivial.  Further, the natural map
	\[
	\iota: \a\times\tilde{\l}\times\l\to\c(\a)
	\]
	is an injective Lie algebra homomorphism, with image $\operatorname{Im}\iota$ an ideal of $\c(\a)$, such that $\operatorname{Im}\iota+\mathfrak{t}=\c(\a)$.
\end{Proposition} 
\begin{proof}
	Since $T'$ is a reduced even root system, the subalgebra $\l$ is a Kac-Moody algebra of finite-type and thus is semisimple.  If $T\neq\emptyset$ then we apply Proposition \ref{prop_about_T}: either $T=\{\pm\alpha\}$ for an odd isotropic root $\alpha$, in which case $\tilde{\l}\cong\s\l(1|1)$, or $T$ is an irreducible GRRS.  The only possibilities for $\tilde{\l}$ in the latter case are then either a basic simple Lie superalgebra or $\s\l(n|n)$ for $n\geq 2$.
	
	Suppose that $\alpha\in T$ is odd and $\beta\in T'$ such that $\alpha+\beta$ is a root; then it is necessarily odd and thus lies in $T$; on the other hand $(\beta,\alpha+\beta)=(\beta,\beta)\neq0$, which would imply that $\beta\in T$, i.e.\ $\beta\notin T'$, a contradiction.   It follows that $[\g_{\alpha},\g_{\beta}]=0$, and thus $[\l,\tilde{\l}]=0$.  
	
	Hence we obtain a natural map $\a\times\tilde{\l}\times\l\to \c(\a)$; the only case it could not be injective is if $\tilde{\l}\cap\a\neq0$; however since $\tilde{\l}\sub\k$ this cannot happen, and we are done.
\end{proof}

\begin{Corollary}\label{centralizater}
	Suppose that $\g$ is a finite dimensional Kac-Moody Lie superalgebra with indecomposable Cartan matrix (i.e.\ we remove the case $\g=\p\s\l(n|n)$).   Then $\c(\a)$ is a product of a reductive Lie algebra with a Kac-Moody Lie superalgebra with an indecomposable Cartan matrix. 
\end{Corollary}

\begin{proof}
	By Lem. 3.1 of \cite{serganova2022splitting}, $\c(\a)$ is the product of an abelian Lie algebra with a symmetrizable Kac-Moody Lie superalgebra.  By Proposition \ref{centralizer_prop}, exactly one factor will have odd roots, so we obtain the result.
\end{proof}

\section{Involutions and the Iwasawa Decomposition}\label{sec_iwasawa}

Let us now assume that $\g$ either is basic simple or is $\g\l(m|n)$ for some $m,n\in\N$, and that $\theta$ is an involution preserving a chosen nondegenerate invariant form on $\g$.  Then in our decomposition $\g=\k\oplus \p$ we have that $\p$ is the $(-1)$-eigenspace of $\theta$.  Recall that on a Lie superalgebra $\g=\g_{\ol{0}}\oplus\g_{\ol{1}}$ there is a canonical involution $\delta\in\Aut(\g)$ defined by $\delta=\id_{\g_{\ol{0}}}\oplus(-\id_{\g_{\ol{1}}})$.  This involution is central in $\Aut(\g)$.

\begin{Lemma}\label{involutions_trivial_even}
	If $\theta\neq\id_{\g},\delta$, then $\p_{\ol{0}}\neq0$.
\end{Lemma}

\begin{proof}
	If $\p_{\ol{0}}=0$, then we have $\g_{\ol{0}}$ is fixed by $\theta$.   Then $\theta$ fixes a Cartan subalgebra $\h\sub\g_{\ol{0}}$, and hence $\theta$ must preserve the root spaces with respect to this Cartan, and so by the assumption that $\theta$ is an involution, it acts by $\pm1$ on each odd root space of $\g$.  Now $\g_{\ol{1}}$ is a $\g_{\ol{0}}$-module, and $\theta$ will be an intertwiner for this module structure.  By the general theory of simple Lie superalgebras (see Chapter 1 of \cite{musson2012lie}), $\g_{\ol{1}}$ either is irreducible or breaks into a sum of two irreducible $\g_{\ol{0}}$-representations $\g_{\ol{1}}',\g_{\ol{1}}''$ such that $[\g_{\ol{1}}',\g_{\ol{1}}'']=\g_{\ol{0}}$ (or $[\g_{\ol{1}}',\g_{\ol{1}}'']$ is a codimension 1 subalgebra of $\g_{\ol{0}}$ in the case of $\g\l(m|n)$).  In the former case, $\theta$ must act by $\pm1$ on $\g_{\ol{1}}$, so that $\theta=\id$ or $\delta$.
	
	In the latter case, let us first assume that $\g\neq\p\s\l(2|2)$ so that $\g_{\ol{1}}'$ and $\g_{\ol{1}}''$ are non-isomorphic $\g_{\ol{0}}$-modules.  If $\theta$ does not act by $\pm1$ on all of $\g_{\ol{1}}$, then WLOG it will act by $(-1)$ on $\g_{\ol{1}}'$ and by $1$ on $\g_{\ol{1}}''$, and thus $[\g_{\ol{1}}',\g_{\ol{1}}'']\sub\p_{\ol{0}}=0$, a contradiction.  
	
	Finally if $\g=\p\s\l(2|2)$, then as is shown in Chpt. 5 of \cite{musson2012lie}, the set of automorphisms of $\g$ that fix $\g_{\ol{0}}$ is isomorphic to $SL_2(\Bbbk)$.  The only order 2 element of $SL_2(\Bbbk)$ is $(-1)$, which corresponds to $\delta$, and so we are done.
\end{proof}

Since we have an involution on $\g_{\ol{0}}$ preserving the nondegenerate form on it, by classical theory (see for instance Sec. 26 of \cite{timashev2011homogeneous}) we may choose a maximal toral subalgebra $\a\sub\p_{\ol{0}}$ that can be extended to a $\theta$-invariant Cartan subalgebra of $\g$, which we will call $\h$.  We obtain a decomposition $\h=\t\oplus\a$, where $\t$ is the fixed subspace of $\theta$.  We again write $\c(\a)$ for the centralizer of $\a$ in $\g$.  Notice that $\a$ is also a maximal toral subalgebra of the $(-1)$-eigenspace of the involution $\delta\circ\theta$.

We already described the structure of $\c(\a)$ as an algebra in Proposition \ref{centralizer_prop}, and in particular we saw that $\c(\a)_{\ol{1}}=\tilde{\l}_{\ol{1}}$.  Now $\theta$ restricts to an automorphism of $\c(\a)$ preserving $\tilde{\l}$, and by classical theory we have $\c(\a)_{\ol{0}}\cap\p=\a$.  Since $\a\cap\tilde{\l}=0$, by Lemma \ref{involutions_trivial_even} either $\theta|_{\tilde{\l}}=\id_{\tilde{\l}}$ or $\theta|_{\tilde{\l}}=\delta_{\tilde{\l}}$.

\begin{Definition}\label{defn_proj}
	For $\lambda\in\h^*$ write $\ol{\lambda}:=(\lambda-\theta\lambda)/2\in\a^*$ for the orthogonal projection of $\lambda$ to $\a^*$ (equivalently the restriction to $\a$), and write $\ol{R}$ for the restriction of roots in $R$ to $\a$ which are nonzero.  We call $\ol{R}\sub\a^*$ the restricted root system of the pair $(\g,\k)$, and elements of $\ol{R}$ we call restricted roots.  
\end{Definition}   
Let $\Z\ol{R}\sub\a^*$ be the $\Z$-module generated by $\ol{R}$, and then choose a group homomorphism $\ol{\phi}:\Z\ol{R}\to\R$ such that $\ol{\phi}(\ol{\alpha})\neq0$ for all $\ol{\alpha}\in \ol{R}$.  Let $\ol{R}^{\pm}=\{\ol{\alpha}\in\ol{R}:\pm\ol{\phi}(\ol{\alpha})>0\}$ so that we obtain a partition of the restricted roots $\ol{R}=\ol{R}^+\sqcup\ol{R}^-$.  We call $\ol{R}^+$ the positive restricted roots, and we call a partition of $\ol{R}$ arising in this way a choice of positive system for $\ol{R}$.  Write 
\[
\n^{\pm}=\bigoplus\limits_{\ol{\alpha}\in\ol{R}^{\pm}}\g_{\ol{\alpha}},
\]
where $\g_{\ol{\alpha}}$ is the weight space of $\ol{\alpha}\in\a^*$ with respect to the adjoint action of $\a$ on $\g$.  We will use $\n=\n^+$ as a shorthand.

\begin{Theorem}\label{main_technical_thm}
	If $\theta|_{\c(\a)_{\ol{1}}}=\id$, then we get an Iwasawa decomposition of $\g$:
	\[
	\g=\k\oplus\a\oplus\n
	\]
\end{Theorem} 

\begin{proof}
	The proof is identical to the classical case.  We see that for $\ol{\alpha}\in \ol{R}$, we have linear isomorphisms $\theta:\g_{\ol{\alpha}}\to\g_{-\ol{\alpha}}$, so that $\g_{\ol{\alpha}}\cap\k=\g_{\ol{\alpha}}\cap\p=0$.  Hence if $y\in\g_{\ol{\alpha}}$ is nonzero and $y=y_0+y_1$ where $y_0\in\k$ and $y_1\in\p$, then $y_0\neq0$ and $y_1\neq0$, and we have $\theta(y)=y_0-y_1$.  From this it is clear that $\k+\a+\n$ contains $\n^-$, and it is also clear that it contains $\h$.  We see $\c(\a)$ is complementary to $\a+\n+\n^{-}$, and by our assumption on $\theta$ we have $\c(\a)\sub\k+\a$, which shows that $\k+\a+\n=\g$.
	
	To show the sum is direct, if we have $x+h+y=0$, where $x\in\k$, $h\in\a$, and $y\in\n$, then applying $[h',\cdot]$ for $h'\in\a$ we find that $[h',y]=-[h',x]\in\p$.  Hence $\theta([h',y])=-[h',y]\in\n$, while $[\theta(h'),\theta(y)]=-[h',\theta(y)]\in\n^-$.  Hence $[h',y]=0$ for all $h'\in\a$ implying $y=0$.  It follows that $x+h=0$, and since $x\in\k$ and $h\in\p$ this implies $x=h=0$, and we are done.
\end{proof}

Before stating the next corollary, we need a couple of definitions.

\begin{Definition}\label{pos_sys_def}
	Let $R$ be a GRRS and let $Q=\Z R\sub\h^*$ be the root lattice.  Given a group homomorphism $\phi:Q\to\R$ such that $\phi(\alpha)\neq0$ for all $\alpha\in R$, we obtain a partition $R=R^+\sqcup R^-$ where $R^{\pm}=\{\alpha\in R:\pm\phi(\alpha)>0\}$.  We call $R^+$ the positive roots of $R$, and any partition of $R$ arising in this way is called a positive system.
\end{Definition}

Positive systems for $R$ are equivalent to choices of Borel subalgebras of the corresponding Lie superalgebra $\g$ which contain $\h$, where the Borel subalgebra is given by $\b=\h\oplus\bigoplus\limits_{\alpha\in R^+}\g_{\alpha}$ (in fact we define Borel subalgebras to be subalgebras arising in this way).

\begin{Definition}
	Let $\theta$ be an involution of $\g$ which admits an Iwasawa decomposition.  We say a choice of positive system for $R$ is an \emph{Iwasawa positive system} if there exists a positive system for $\ol{R}$ which is compatible with it.  Here, if $\ol{\phi}:\ol{R}\to\R$ and $\phi:R\to\R$ are homomorphisms determining positive systems for $\ol{R}$ and $R$ respectively, we say $\phi$ is compatible with $\ol{\phi}$ if $\phi(\alpha)>0$ whenever both $\ol{\alpha}\neq0$ and $\ol{\phi}(\ol{\alpha})\neq0$.  If $R$ is an Iwasawa positive system, we call the corresponding Borel subalgebra an Iwasawa Borel subalgebra of $\g$.
\end{Definition}

\begin{Corollary}\label{main_cor}
	If $\theta$ is an involution on a basic simple Lie superalgebra or $\g\l(m|n)$ such that $\theta$ preserves the nondegenerate invariant form, then the following are true:
	\begin{enumerate}
		\item either $\theta$ or $\delta\circ\theta$ admits an Iwasawa decomposition;
		\item an Iwasawa Borel subalgebra of $\g$ corresponding to $\theta$ exists, and it is complementary to the fixed points of $\theta$ if $\theta$ admits an Iwasawa decomposition.
	\end{enumerate} 
\end{Corollary}

\begin{proof}
	By Proposition \ref{centralizer_prop}, either $\theta$ or $\delta\circ\theta$ satisfies the hypothesis of Theorem \ref{main_technical_thm}. If $\g=\p\s\l(2|2)$ we reference the classification of involutions in \cite{serganova1983classification}. 
	
	To construct an Iwasawa Borel subalgebra, we construct an Iwasawa positive system.  Let $\ol{\phi}:\Z\ol{R}\to\R$ be a group homomorphism determining a positive system for $\ol{R}$.  Split the natural surjection of free abelian groups $\Z R\to\Z\ol{R}$ so that $\Z R\cong\Z\ol{R}\oplus K$.  Then construct $\phi:\Z R\to \R$ which is an extension of $\ol{\phi}$ with respect to the inclusion $\Z\ol{R}\to\Z R$, such that both $\phi(\alpha)\neq0$ for all $\alpha\in R$ and $\phi(\alpha)>0$ whenever $\ol{\phi}(\ol{\alpha})>0$ for $\alpha\in R$.  Then the Iwasawa Borel subalgebra $\b=\h\oplus\bigoplus\limits_{\phi(\alpha)>0}\g_{\alpha}$ contains $\a\oplus\n$ and thus is complementary to $\k$ by the Iwasawa decomposition.
\end{proof} 

\begin{Proposition}\label{cent_gen_simple_roots}
	Let $\theta$ be an involution as in Corollary \ref{main_cor} and suppose that $\b$ is an Iwasawa Borel subalgebra of $\g$.  Then the simple roots of $\b$ that are fixed by $\theta$ generate all fixed roots of $\theta$.  In particular, $\c(\a)$ is generated by $\h\sqcup\{e_{\gamma},e_{-\gamma}\}_{\gamma\in I}$, where $I$ is the set of positive simple roots fixed by $\theta$.
\end{Proposition}

\begin{proof}
	If $\beta$ is a positive root then we may write
	\[
	\beta=\sum\limits_{\alpha\notin I} c_\alpha\alpha+\sum\limits_{\gamma\in I} d_\gamma\gamma
	\]
	where the first sum is over simple roots $\alpha$ not fixed by $\theta$, and $c_\alpha,d_{\gamma}\in\Z_{\geq0}$.  If $\theta\beta=\beta$ then we obtain that
	\[
	\beta=\sum\limits_{\alpha\notin I} c_\alpha\theta\alpha+\sum\limits_{\gamma\in I} d_\gamma\gamma.
	\]
	But $\theta\alpha$ is a negative root for $\alpha\notin I$, and thus $c_{\alpha}=0$.
\end{proof}

  We give a list of supersymmetric pairs for the Lie superalgebra $\g\l(m|n)$ and the basic simple Lie superalgebras not of type $A$.  We first describe their generalized root systems explicitly.  
\begin{itemize}
	\item $\g=\g\l(m|n)$: $V=\Bbbk\langle\epsilon_1,\dots,\epsilon_m,\delta_1,\dots,\delta_n\rangle$, $(\epsilon_i,\epsilon_j)=-(\delta_i,\delta_j)=\delta_{ij}$, $(\epsilon_i,\delta_j)=0$.  The even roots are 
	\[
	R_{ev}=\{\epsilon_i-\epsilon_j:i\neq j\}\sqcup\{\delta_i-\delta_j:i\neq j\}
	\]
	and the odd roots are 
	\[
	R_{odd}=\{\pm(\epsilon_i-\delta_j)\}.
	\]

	\item $\g=\o\s\p(2m|2n)$: $V=\Bbbk\langle\epsilon_1,\dots,\epsilon_m,\delta_1,\dots,\delta_n\rangle$, $(\epsilon_i,\epsilon_j)=-(\delta_i,\delta_j)=\delta_{ij}$, $(\epsilon_i,\delta_j)=0$.  The even roots are 
	\[
	R_{ev}=\{\pm\epsilon_i\pm\epsilon_j:i\neq j\}\sqcup\{\pm\delta_i\pm\delta_j:i\neq j\}\sqcup\{\pm2\delta_i\}
	\] 
	and the odd roots are 
	\[
	R_{odd}=\{\pm\epsilon_i\pm\delta_j\}
	\]
	
	\item $\g=\o\s\p(2m+1|2n)$: $V=\Bbbk\langle\epsilon_1,\dots,\epsilon_m,\delta_1,\dots,\delta_n\rangle$, $(\epsilon_i,\epsilon_j)=-(\delta_i,\delta_j)=\delta_{ij}$, $(\epsilon_i,\delta_j)=0$.  The even roots are 
	\[
	R_{ev}=\{\pm\epsilon_i\pm\epsilon_j:i\neq j\}\sqcup\{\pm\epsilon_i\}\sqcup\{\pm\delta_i\pm\delta_j:i\neq j\}\sqcup\{\pm2\delta_i\}
	\]
	and the odd roots are 
	\[
	R_{odd}=\{\pm\epsilon_i\pm\delta_j\}\sqcup\{\pm\delta_i\}.
	\]
	
	\item $\g=D(2,1;a)$: $V=\Bbbk\langle\epsilon,\delta,\gamma\rangle$, $(\epsilon,\epsilon)=1$, $(\delta,\delta)=a$, $(\gamma,\gamma)=-a-1$, and $(\epsilon,\delta)=(\epsilon,\gamma)=(\delta,\gamma)=0$.  The even roots are 
	\[
	R_{ev}=\{\pm2\epsilon,\pm 2\delta,\pm2\gamma\}, 
	\]
	and the odd roots are
	\[
	R_{odd}=\{\pm\epsilon\pm\delta\pm\gamma\}.
	\]
	
	\item $\g=\a\b(1|3)$, root system is $AB(1|3)$: $V=\Bbbk\langle\delta,\epsilon_1,\epsilon_2,\epsilon_3\rangle$, $(\delta,\delta)=-1$, $(\epsilon_i,\epsilon_j)=\delta_{ij}/3$.  The even roots are 
	\[
	R_{ev}=\{\pm\delta\}\sqcup\{\pm\epsilon_i,\pm\epsilon_i\pm\epsilon_j:i\neq j\}
	\]
	and the odd roots are 
	\[
	R_{odd}=\{\frac{1}{2}(\pm\delta\pm\epsilon_1\pm\epsilon_2\pm\epsilon_3)\}.
	\]
	
	\item $\g=\a\g(1|2)$: $V=\Bbbk\langle\delta,\epsilon_1,\epsilon_2,\epsilon_3\rangle$ with the relation $\epsilon_1+\epsilon_2+\epsilon_3=0$, and inner product $(\epsilon_i,\epsilon_i)=-2(\epsilon_i,\epsilon_j)=-(\delta,\delta)=2$, where $i\neq j$.  Then the even roots are 
	\[
	R_{ev}=\{\pm\epsilon_i,\epsilon_i-\epsilon_j:i\neq j\}\sqcup\{\pm2\delta\}
	\]
	and odd roots 
	\[
	R_{odd}=\{\pm\delta\}\sqcup\{\pm\delta\pm\epsilon_j\}.
	\]
\end{itemize}  

We will now give a list of supersymmetric pairs for each of the superalgebras in the above list.  For superalgebras not of type $A$, we will give all supersymmetric pairs up to conjugacy of the corresponding involution.  For those of type $A$ we will only describe two families of pairs for $\g\l(m|n)$, since these are the most prominent in the literature and are exactly those which lift to the supergroup $GL(m|n)$.  Further, any other supersymmetric pair for $\g\l(m|n)$ determined by an involution $\theta$ is conjugate to one of these two families up to its action on the center of $\g\l(m|n)$.  The supersymmetric pairs for $\s\l(m|n)$ with $m\neq n$ and $\p\s\l(n|n)$ coming from involutions preserving an invariant form come from these two families for $\g\l(m|n)$; for precise details, see \cite{serganova1983classification}.

For a proof of the statement that these are all such supersymmetric pairs when $\g\neq\a\g(1|2)$ or $\a\b(1|3)$, we refer to Serganova's classification in \cite{serganova1983classification}.   The cases for $\a\g(1|2)$ and $\a\b(1|3)$ were communicated to the author by Serganova, and are written in the appendix.

In each case of the below table we describe the action of the involution on basis elements of the GRRS, where we omit any basis elements that are fixed by the involution. For cases (1) and (3) we are giving the GRRS automorphism when $r\leq m/2$ and $s\leq n/2$.

\renewcommand{\arraystretch}{2}

\begin{center}
\begin{tabular}{|c|c|c|}
	\hline 
	Supersymmetric Pair & Iwasawa Decomposition? & GRS Automorphism \\
	\hline
	\makecell{$(\g\l(m|n)$,\\ $\g\l(r|s)\times\g\l(m-r|n-s))$} & Iff $(m-2r)(n-2s)\geq0$   & \makecell{$\epsilon_i\leftrightarrow\epsilon_{m-i+1}, 1\leq i\leq r$,\\ $\delta_j\leftrightarrow\delta_{n-j+1}, 1\leq j\leq s$} \\
	\hline
	$(\g\l(m|2n),\o\s\p(m|2n))$ & Yes &  $\epsilon_i\leftrightarrow-\epsilon_i, \delta_i\leftrightarrow-\delta_{2n-i+1}$ \\
	\hline
	\makecell{$(\o\s\p(m|2n)$,\\$\o\s\p(r|2s)\times$ \\ $\o\s\p(m-r,2n-2s))$}& Iff $(m-2r)(n-2s)\geq0$  & \makecell{$\epsilon_i\leftrightarrow-\epsilon_i, 1\leq i\leq r$\\ $\delta_i\leftrightarrow\delta_{n-i+1}, 1\leq i\leq s$}\\
	\hline
	$(\o\s\p(2m|2n),\g\l(m|n))$  & Yes & $\delta_i\leftrightarrow-\delta_i, \epsilon_i\leftrightarrow-\epsilon_{m-i+1}$ \\
	\hline
	$(D(2,1;a),\o\s\p(2|2)\times\s\o(2))$ &  Yes & $\epsilon\leftrightarrow-\epsilon,\delta\leftrightarrow-\delta$  \\
	\hline
	 $(\a\b(1|3),\s\l(1|4))$ & Yes & $\epsilon_1\leftrightarrow	-\epsilon_1,\delta\leftrightarrow-\delta$\\
	\hline
	$(\a\b(1|3),\g\o\s\p(2|4))$ & Yes & \makecell{$\epsilon_1\leftrightarrow-\epsilon_1,\epsilon_2\leftrightarrow-\epsilon_2,$\\ $\delta\leftrightarrow-\delta$}\\
	\hline
	$(\a\b(1|3),D(2,1;2)\times \s\l(2))$ & Yes & $\epsilon_i\leftrightarrow-\epsilon_i$\text{ for all }$i$\\
	\hline
	$(\a\g(1|2),D(2,1;3))$ & Yes &  $\epsilon_i\leftrightarrow-\epsilon_i$\text{ for all }$i$\\
	\hline
	$(\a\g(1|2),\o\s\p(3|2)\times\s\l_2)$ & No &  $\epsilon_i\leftrightarrow-\epsilon_i$\text{ for all }$i$\\
	\hline
\end{tabular}
\end{center}

\

\

Note that $\o\s\p(1|2)$ does not admit a nontrivial involution preserving the form by Lemma \ref{no_odd_negation}, and thus by Lemma \ref{killing_form_pres} has no nontrivial involutions.  Further Lemma \ref{no_odd_negation} also implies there is never an involution that acts by (-1) on a Cartan subalgebra and preserves the form.  This may seem surprising given the existence of the Chevalley involution for reductive Lie algebras.  The following remark seeks to contextualize this.
\begin{remark}	
	A complex Kac-Moody Lie algebra $\g$ always admits a nontrivial involution $\omega$, the Chevalley involution, that acts by $(-1)$ on a Cartan subspace (see \cite{kac1990infinite} Chapter 1).  If one modifies this involution to make it complex antilinear as in Chapter 2 of \cite{kac1990infinite}, one can construct a Cartan involution of $\g$, i.e.\ an involution whose fixed points are a compact real form of $\g$.  For finite type complex Kac-Moody algebras one can use Cartan involutions to set up a bijection between real forms of $\g$ and complex linear involutions of $\g$, as originally shown by Cartan.  
	
	For complex Kac-Moody Lie superalgebras the natural generalization of the Chevalley involution which we write as $\tilde{\omega}$, is of order 4.  In fact $\tilde{\omega}^2=\delta$, so it is of order 2 on $\g_{\ol{0}}$ and order 4 on $\g_{\ol{1}}$.  Write $\Aut_{2,4}(\g)$ for the complex linear automorphisms $\theta$ of $\g$ which are order $2$ on $\g_{\ol{0}}$ and order $4$ on $\g_{\ol{1}}$.  Then if $\g$ a finite-dimensional contragredient Lie superalgebra then there is a bijection between the real forms of $\g$ and $\Aut_{2,4}(\g)$ as shown in \cite{chuah2013cartan}.
\end{remark}

\begin{remark}\label{counterexample}
	There are other supersymmetric pairs for the algebras we consider that are often studied but which do not appear in the list above -- namely $(\g\l(n|n),\p(n))$ and $(\g\l(n|n),\q(n))$.  However these are exactly the cases when the involution does not preserve an invariant form, which can be seen from the fact that neither $\p(n)$ nor $\q(n)$ are basic.  For the pair $(\g\l(n|n),\q(n))$ the Iwasawa decomposition does hold as the Cartan subspace in that case contains a regular semisimple element.
	
	However Proposition \ref{centralizer_prop} and in particular Corollary \ref{main_cor} fail for the pair $(\g\l(n|n),\p(n))$.  We will show this now, and it demonstrates the necessity of the automorphism to preserve the form.  The involution in this case, which we call $\theta$, is given explicitly by
	\[
	\begin{bmatrix}
		W & X\\ Y & Z
	\end{bmatrix}\mapsto \begin{bmatrix}
		-Z^t & X^t\\-Y^t & -W^t
	\end{bmatrix}
	\]
	Thus a Cartan subspace is given by
	\[
	\a=\left\{\begin{bmatrix}
		D & 0\\ 0 & D
	\end{bmatrix}:D\text{ is diagonal}\right\}.
	\]
	Hence
	\[
	\c(\a)=\left\{\begin{bmatrix}
		D & D'\\ D' & D
	\end{bmatrix}:D,D'\text{ are diagonal}\right\}\cong\s\l(1|1)\times\cdots\s\l(1|1).
	\]
	So Proposition \ref{centralizer_prop} fails.  Further we see that $\theta|_{\c(\a)_{\ol{1}}}\neq \pm\id_{\c(\a)_{\ol{1}}}$, so Corollary \ref{main_cor} fails too.  In particular $(\g\l(n|n),\p(n))$ does not admit an Iwasawa decomposition.
	
	However despite the failure of having an Iwasawa decomposition, $\p(n)$ is still a spherical subalgebra of $\g\l(n|n)$, i.e.\ there is a complementary Borel subalgebra to $\p(n)$ in $\g\l(n|n)$.  In particular the Borel subalgebra with simple roots $\delta_1-\epsilon_1,\epsilon_1-\delta_2,\delta_2-\epsilon_2,\dots,\epsilon_{n-1}-\delta_n,\delta_n-\epsilon_n$ is complementary to $\p(n)$.  In fact, this is the only Borel subalgebra with this property up to conjugacy, i.e.\ up to inner automorphisms.  
	
	Indeed, if $\b$ is a such a Borel subalgebra then we may decompose it according to its $\Z$-grading as $\b=\b_{-1}\oplus\b_0\oplus\b_1$, induced by the $\Z$-grading on $\g\l(m|n)$ (coming from a grading operator).  Then by dimension reasons we must have $\dim\b_{-1}=n(n+1)/2$ and $\dim\b_{1}=n(n-1)/2$.  Using the indexing of conjugacy classes of Borel subalgebras of $\g\l(m|n)$ by $\epsilon\delta$-sequences as explained in Sec.\ 1.3 of \cite{cheng2012dualities}, one can see that there is a unique conjugacy class of Borel subalgebras with these dimensions for $\b_{\pm1}$, giving us uniqueness.
\end{remark}

For the superalgebras we consider, a choice of positive system is equivalent to a choice of simple roots in the GRRS, just as with even root systems. 

A choice of simple roots can be encoded in a Dynkin-Kac diagram, and one obtains a bijection between Dynkin-Kac diagrams and choices of simple roots up to Weyl group symmetries for a given superalgebra (see \cite{kac1977lie}).  Just as in the classical case, if one chooses an Iwasawa positive system, one can construct a \emph{Satake diagram} from it using the results of the following lemma, which is proven exactly as in \cite{satake1960representations}. For this subsection we only consider one of the supersymmetric pairs in our table above,  so that simple roots are linearly independent.
\begin{Lemma}\label{diagram_auto}
	Let $\Pi$ be the set of simple roots of an Iwasawa positive system.  If $\alpha$ is a simple root such that $\theta\alpha\neq\alpha$, then 
	\[
	-\theta\alpha=\alpha'+\sum_{\gamma\in I}d_\gamma\gamma
	\]
	where $\alpha'$ is a simple root and $I\sub\Pi$ is the set of simple roots fixed by $\theta$.  The correspondence $\alpha\mapsto\alpha'$ defines an permutation of order 1 or 2 on $\Pi\setminus I$.  In particular, for distinct simple roots $\alpha,\beta$, we have $\ol{\alpha}=\ol{\beta}$ (see Definition \ref{defn_proj} for the notation $\ol{\alpha},\ol{\beta}$) if and only if $\beta=\alpha'$.
\end{Lemma}
\begin{proof}
	Write $\{\alpha_i\}_{i}$ for the set of simple roots not fixed by $\theta$.  Then $-\theta\alpha_i$ is a positive root for all $i$, and thus we may write
	\[
	-\theta\alpha_i=\sum\limits_{j}c_{ij}\alpha_j+\sum\limits_{\gamma\in I}d^{i}_{\gamma}\gamma
	\]
	for some $d^{i}_\gamma\in\Z_{\geq0}$, where $C=(c_{ij})$ is square and has nonnegative integer entries.  Applying $(-\theta)$ to this equation once again, we obtain that
	\[
	\alpha_i=\sum\limits_{j,k}c_{ij}c_{jk}\alpha_k+\sum r^i_\gamma\gamma
	\]
	for some $r^i_{\gamma}\in\Z$.  Since $\alpha_i$ is simple, this forces $C^2$ to be the identity matrix, which implies that $C$ is in fact a permutation matrix.  This permutation matrix defines our permutation of $\Pi\setminus I$.	
	
	For the last statement, if $\ol{\alpha}=\ol{\beta}$, then $\alpha-\theta\alpha=\beta-\theta\beta$, so there exists $\gamma_\alpha,\gamma_\beta$ in the span of fixed simple roots such that
	\[
	\alpha+\alpha'+\gamma_\alpha=\beta+\beta'+\gamma_\beta.
	\]
	By linear independence of our base, we must have that $\{\alpha,\alpha'\}=\{\beta,\beta'\}$, so we are done.
\end{proof}
Using the above result, we may construct a Satake diagram from $(\g,\k)$ as follows: choosing an Iwasawa positive system, we get a Dynkin-Kac diagram for $\g$.  Now draw an arrow between two distinct simple roots if they are related by the involution constructed in Lemma \ref{diagram_auto}.   Finally, draw a solid black line over a node if the corresponding simple root $\alpha$ is fixed by $\theta$.  Classically one would color the node black, but unfortunately Dynkin-Kac diagrams may already have black nodes as they represent non-isotropic odd simple roots.  

We call the result a Satake diagram for the corresponding supersymmetric pair.  Note that it is not unique-- Proposition \ref{Iwasawa_systems} shows that it is determined exactly up to choices of positive systems for $\ol{R}$ and $S$ (see Section 6 for more on the structure of $\ol{R}$).  Others have given examples of such diagrams, such as in \cite{pati1998satake}.  In that paper nodes are drawn black if the corresponding simple root is fixed by $\theta$.  

Before we state the proposition, we define a positive system of $S$ to be a choice of positive and negative roots in $S$ arising from a group homomorphism $\psi:\Z S\to\R$ such that $\psi(\gamma)\neq0$ for all $\gamma\in S$, as in Definition \ref{pos_sys_def} (recall $S$ might not be a GRRS). 
\begin{Proposition}\label{Iwasawa_systems}
	There is a natural bijection between Iwasawa positive systems and choices of positive systems for $\ol{R}$ and $S$.
\end{Proposition}
\begin{proof}
	The simple roots of any positive root system form a $\Z$-basis of $Q:=\Z R$, the root lattice.  Thus by Lemma \ref{cent_gen_simple_roots} we have that $\Z S$ splits off from $Q$, so we can write $Q=\Z S\oplus Q'$.  Write $\pi:Q\to\Z\ol{R}$ for the canonical projection, and observe that $\Z S\sub\ker\pi$.  Therefore the restricted map $Q'\to\Z\ol{R}$ is surjective, so we may split it and write $Q'=\Z\ol{R}\oplus Q''$, so that $Q'=\Z S\oplus\Z\ol{R}\oplus Q''$.
	
	Now let $\phi:Q\to\R$ be a group homomorphism determining an Iwasawa positive system coming from $\ol{\phi}:\Z\ol{R}\to\R$ as in Corollary \ref{main_cor}.  Write $\psi:\Z S\to\R$ for the restriction of $\phi$ to $\Z S$.  Then since $\psi(\gamma)\neq0$ for all $\gamma\in S$, $\psi$ determines a positive system for $S$. Thus the Iwasawa positive system gives rise to positive systems of $\ol{R}$ and $S$ respectively from $\ol{\phi}$ and $\psi$.
	
	Conversely, given positive systems of $\ol{R}$ and $S$ coming from group homomorphisms $\ol{\phi}:\Z\ol{R}\to\R$ and $\psi:\Z S\to\R$, the map $\phi:\Z R\to\R$ defined by $\phi=\epsilon\psi\oplus\ol{\phi}\oplus0:\Z S\oplus\Z\ol{R}\oplus Q''\to\R$ determines an Iwasawa positive system, where $\epsilon=\frac{\min_{\ol{R}^+}(\ol{\phi})}{2\max_{S}(\psi)}$. The described correspondences are seen to be bijective, and thus we are done.
\end{proof}
\section{Restricted Root Systems}

Consider one of the supersymmetric pairs $(\g,\k)$ from the table of Section \ref{sec_iwasawa} which admits an Iwasawa decomposition.  Write $\theta$ for the involution, and by abuse of notation also write $\theta$ for the induced involution on the GRRS $R\sub\h^*$ coming from a Cartan subalgebra $\h$ containing a Cartan subspace $\a$.  Continue writing $Q=\Z R\sub\h^*$ for the root lattice, $S\sub R$ for the roots fixed by $\theta$, and $\ol{R}$ for the restricted roots.  We make a few notes about differences between the super case and the purely even case.

For an even symmetric pair there are often roots $\alpha$ for which $\theta(\alpha)=-\alpha$.  In the super case this cannot hold for odd roots.
\begin{Lemma}\label{no_odd_negation}
	If $\alpha$ is an odd root, then $\theta(\alpha)\neq-\alpha$.
\end{Lemma}
\begin{proof}
	Suppose $\alpha$ is odd and satisfies $\theta(\alpha)=-\alpha$.  Write $h_{\alpha}\in\h$ for the coroot of $\alpha$, i.e.\ $h_{\alpha}$ satisfies $(h_{\alpha},-)=\alpha$ as an element of $\h^*$. Then we may assume $\theta e_{\alpha}=e_{-\alpha}$ and $\theta e_{-\alpha}=e_{\alpha}$ where $e_{\alpha}\in\g_{\alpha}$, $e_{-\alpha}\in\g_{-\alpha}$ are nonzero and satisfy $[e_{\alpha},e_{-\alpha}]=[e_{-\alpha},e_{\alpha}]=h_{\alpha}$.  But then 
	\[
	\theta h_{\alpha}=\theta[e_{\alpha},e_{-\alpha}]=[\theta e_{\alpha},\theta e_{-\alpha}]=[e_{-\alpha},e_{\alpha}]=h_{\alpha}.
	\]
	However the action of $\theta$ on $\h^*$ is dual to the action of $\theta$ on $\h$, so since $\alpha$ and $h_{\alpha}$ are dual to one another we must have $\theta h_{\alpha}=-h_{\alpha}$, a contradiction.
\end{proof}

Another proof of the above result can be given by using that $(-,\theta(-))$ defines a nondegenerate symplectic form on $(\g_{\ol{\alpha}})_{\ol{1}}$ for a restricted root $\ol{\alpha}\in\ol{R}$.  Thus $\dim\g_{\ol{\alpha}}$ must be even, so the GRRS involution $(-\id)\circ\theta$ cannot fix any odd roots.

The following lemma is well-known from the even case, and is proven in \cite{araki1962root}.  
\begin{Lemma}\label{even_roots_nice}
	If $\alpha$ is an even root, then $\theta\alpha+\alpha$ is not a root.
\end{Lemma}

However the corresponding statement for odd roots is false in many cases, for instance, for the pair $(\g\l(m|2n),\o\s\p(m|2n))$, $\theta\alpha+\alpha$ is always a root if $\alpha$ is odd.   This property sometimes hold and sometimes fails for other pairs.

Classically, $\ol{R}$ defines a (potentially non-reduced) root system in $\a^*$, the restricted root system of the symmetric pair.  Each restricted root $\ol{\alpha}$ has a positive integer multiplicity attached to it given by $m_{\ol{\alpha}}:=\dim\g_{\ol{\alpha}}$.  The data of the restricted root system with multiplicities completely determines the corresponding symmetric pair.

In the super case it is less clear what type of object the restricted root system is.  Even and odd roots can restrict to the same element of $\a^*$, so the natural replacement of the multiplicity of a restricted root is (a multiple of) the superdimension of the corresponding weight space.  In many cases the object obtained behaves like a GRRS from a combinatorial perspective, however the bilinear form is deformed.  We discuss this situation at the end of of this section, but first state what can be proven in general.

Set $\ol{R}_{re}=\{\ol{\alpha}:\alpha\in R_{re}, \ol{\alpha}\neq0\}\sub\a^*$, $\ol{R}_{im}=\ol{R}\setminus\ol{R}_{re}$.
\begin{Proposition}\label{res_root_sys_props}
	The set $\ol{R}\sub\a^*$ with the restricted bilinear form satisfies the following properties:
	\begin{enumerate}
		\item $\span\ol{R}=\a^*$ and $\ol{R}=-\ol{R}$;
		\item The form is nondegenerate;
		\item Given $\ol{\alpha}\in\ol{R}_{re}$, we have $k_{\ol{\alpha},\ol{\beta}}:=2\frac{(\ol{\alpha},\ol{\beta})}{(\ol{\alpha},\ol{\alpha})}\in\Z$ and $r_{\ol{\alpha}}(\ol{\beta})=\ol{\beta}-k_{\ol{\alpha},\ol{\beta}}\ol{\alpha}\in\ol{R}$.
		\item Given $\ol{\alpha}\in\ol{R}_{im}$, $\ol{\beta}\in\ol{R}$ with $\ol{\beta}\neq\pm\ol{\alpha}$, if $(\ol{\alpha},\ol{\beta})\neq0$ then at least one of $\ol{\beta}\pm\ol{\alpha}\in\ol{R}$.
	\end{enumerate}
	Further, $\ol{R}_{re}\sub \operatorname{span}(\ol{R}_{re})$ is an even (potentially non-reduced) root system and $\ol{R}_{im}$ is invariant under its Weyl group.
\end{Proposition}
\begin{proof}
	Property (1) is obvious, and (2) follows from the fact that we are only considering Lie superalgebras with nondegenerate invariant forms and our involution preserves the form.  The statement (3) is proven just as in the classical case.  For (4), since $(\ol{\alpha},\ol{\beta})\neq0$,  either $(\alpha,\beta)\neq0$ or $(-\theta\alpha,\beta)\neq0$ so either $\beta\pm\alpha$ or $\beta\pm(-\theta\alpha)$ is a root. Restricting to $\a$ gives the desired statement.
	
	That $\ol{R}_{re}$ is a root system is classical (see for instance Chapter 26 of \cite{timashev2011homogeneous}), and it is easy to see that $\ol{R}_{im}$ is Weyl group invariant. 
\end{proof}
\begin{remark}
	Although we use the notation $\ol{R}_{im}$, it is not true in general that $(\ol{\alpha},\ol{\alpha})=0$ for $\ol{\alpha}\in\ol{R}_{im}$.  This is a prominent difference between restricted root systems and GRRSs.
\end{remark}

Using Proposition \ref{res_root_sys_props} we may now $\ol{R}_{re}$ into a union of irreducible real root systems, $\ol{R}_{re}=\ol{R}_{re}^1\sqcup\cdots\ol{R}_{re}^k$.  Since $R$ was irreducible we know that $k\leq 3$ by Proposition \ref{components_prop}.  We may decompose $\a^*$ as $\a^*=U_0\oplus U_1\oplus\cdots\oplus U_k$, where $U_i=\span(\ol{R}_{re}^i)$, and we set $U_0=(\sum_{i\geq 1}U_i)^{\perp}$.  Write $p_i:\a^*\to U_i$ for the projection maps.  The following result is obvious.
\begin{Lemma}\label{restricted_real_comps}
	A real component $\ol{R}^i_{re}$ of $\ol{R}$ either is gotten by
	\begin{itemize}
		\item[(1)] the restriction of nonisotropic roots in a real component of $R_{re}$ preserved by $\theta$; or
		\item[(2)] is obtained as a diagonal subspace of two isomorphic real components of $R$ that are identified by $\theta$. 
	\end{itemize} 
\end{Lemma}
From Lemma \ref{restricted_real_comps} we can prove:
\begin{Proposition}\label{small_orbits_res_roots}
	For each $i>0$, $q_i(\ol{R}_{im})\setminus\{0\}$ is a union of small $W_i$-orbits. 
\end{Proposition}
\begin{proof}
	Let $\ol{\alpha},\ol{\beta}\in\ol{R}_{im}$ such that $q_i(\ol{\alpha}),q_{i}(\ol{\beta})\neq0$ and they lie in the same $W_i$-orbit.  Let $\alpha,\beta\in R$ be lifts of $\ol{\alpha}$ and $\ol{\beta}$.
	
	If $\ol{R}_{re}^i$ falls into the second case of Lemma \ref{restricted_real_comps}, then if we write $p$ for the projection from $\h^*$ onto one of the real components being folded into $\ol{R}_{re}^i$ then $p\alpha$ and $p\beta$ must be conjugate under the Weyl group for that real component too, so we can apply Lemma \ref{projection_non_isotropic}.
	
	Suppose on the other hand that $\ol{R}_{re}^i$ falls into the first case of Lemma \ref{restricted_real_comps}.  Write $p$ for the projection from $\h^*$ onto the corresponding real component giving $\ol{R}_{re}^i$.  Then if $p\alpha$ and $p\beta$ are conjugate under the Weyl group we can apply Lemma \ref{projection_non_isotropic}.  If they are not conjugate under the Weyl group, $R$ must have two imaginary components (see Definition \ref{sec_defs}).  If $\theta$ preserves the imaginary components, then $\alpha$ and $-\theta\beta$ will lie in the same imaginary component and project to $\ol{\alpha},\ol{\beta}$ still, so we are done.  If $\theta$ permutes the imaginary components, then the supersymmetric pair either is $(\g\l(m|2n),\o\s\p(m|2n))$ or $(\o\s\p(2|2n),\o\s\p(1|2n-2r),\o\s\p(1|2r))$.  In the first case $\ol{R}_{re}=A_{m-1}\sqcup A_{n-1}$ so that $\alpha$ and $\beta$ cannot be in distinct imaginary components of $R$, and in the second case $p(R_{im})$ is a single small Weyl group orbit anyway. 
	
\end{proof}
\begin{remark}
	It may be interesting to classify all root systems satisfying the above properties.  That is we consider a complex inner product space $V$ with a finite set $R\sub V$ partitioned into real and imaginary roots $R=R_{re}\sqcup R_{im}$ such that all the properties of Proposition \ref{res_root_sys_props} and Lemma \ref{small_orbits_res_roots} hold.  We will call such objects restricted generalized root systems (RGRSs).  We can ask what all (irreducible) RGRSs are.
	
	Amongst them we would have all deformed weak generalized root systems (WGRSs) as defined in below.  However there would be more examples. One interesting case (communicated to the author by Serganova) comes from the supersymmetric pair $(\a\b(1|3),D(2,1;2))$ where the restricted root system has $R_{re}=B_3$ and $R_{im}=W\omega_{3}$, where $\omega_3$ is the fundamental weight giving the spinor representation of $\s\o(7)$.  
	
	Another exotic example would be $V=\Bbbk^4$, $R_{re}=A_1\sqcup A_1\sqcup A_1\sqcup A_1$ and $R_{im}=W(\omega_1^{(1)}+\omega_1^{(2)}+\omega_1^{(3)}+\omega_1^{(4)})$ where the inner product on each $A_1$ is the same.  This case has four real components which cannot happen for a GRRS.  However one can show that an irreducible RGRS can have at most four components.
\end{remark}

In the case when $\ol{R}_{re}$ has more than one component, it turns out that the restricted root system is a deformed GRS, as introduced in \cite{sergeev2004deformed}.  There, they introduce generalized root systems as more a general object than in \cite{serganova1996generalizations} by relaxing condition (4) in Definition \ref{GRRS_def} to
\begin{itemize}
	\item[(4')] If $\alpha,\beta\in R$ and $(\alpha,\alpha)=0$, then if $(\alpha,\beta)\neq0$ at least one of $\beta\pm\alpha\in R$.  
\end{itemize}
It is also assumed that the inner product is nondegenerate.  It is shown in \cite{serganova1996generalizations} that in a GRRS only one of $\beta\pm\alpha$ can be in $R$.  Following \cite{gorelik2017generalized}, we will call the notion of GRS in the sense of \cite{sergeev2004deformed} a weak GRS (WGRS).  Serganova classified all WGRSs in Section 7 of \cite{serganova1996generalizations}; there are two cases that do not appear in the classification of GRSs:
\begin{itemize}
	\item $C(m,n)$, $m,n\geq 1$: $R_{re}^1=C_m$, $R_{re}^2=C_n$, $R_{im}=W(\omega_1^{(1)}+\omega_1^{(2)})$
	\item $BC(m,n)$, $m,n\geq 1$: $R_{re}^1=BC_m$, $R_{re}^2=C_n$, $R_{im}=W(\omega_1^{(1)}+\omega_1^{(2)})$.
\end{itemize}
Sergeev and Veselov define a deformed WGRS as the data of a WGRS with a deformed inner product determined by a nonzero parameter $k\in \Bbbk^\times$, along with Weyl-group invariant multiplicities $m(\alpha)\in \Bbbk$ for each root $\alpha\in R$.  These multiplicities are required to satisfy certain polynomial relations and that $m(\alpha)=1$ for an isotropic (with respect to the non-deformed bilinear form) root $\alpha$.  

We now explain when and how we can realize $\ol{R}$ as a deformed WGRS.  For each of the supersymmetric pairs we consider where $\ol{R}_{re}$ has more than one component the deformation parameter $k$ is determined by the restriction of the form.  In this case $\ol{R}_{im}\neq\emptyset$, and the multiplicity of every $\ol{\alpha}\in\ol{R}_{im}$ is $-\ell$ for some positive integer $\ell$.  We define the multiplicities of a restricted root $\ol{\alpha}\in R$ to be $m(\ol{\alpha})=-\frac{1}{\ell}\operatorname{sdim}\g_{\ol{\alpha}}$.  Then we claim that we obtain a deformed WGRS in this way.  This can be checked case by case, and we do this in the table below.  Note that this fact has been known to several researchers for some time (most of whom knew before the author).  We give this information here for the benefit of the reader.

In the table below we list, for each supersymmetric pair we consider in which $\ol{R}_{re}$ has more than one component, the corresponding Sergeev-Veselov deformation parameters.  

\

\begin{center}
\begin{tabular}{|c|c|c|c|c|c|}
	\hline
	Supersymmetric Pair & $k$ & $p$ & $q$ & $r$ & $s$ \\
	\hline
	\makecell{$(\g\l(m|n)$,\\ $\g\l(r|s)\times\g\l(m-r|n-s))$} & $-1$ & \makecell{$(n-m)+$\\ $2(r-s)$} & $-\frac{1}{2}$ & \makecell{$(m-n)+$\\ $2(s-r)$} & $-\frac{1}{2}$   \\
	\hline
	$(\g\l(m|2n),\o\s\p(m|2n))$ & $-\frac{1}{2}$ & 0 & 0 & 0 & 0 \\
	\hline
	\makecell{$(\o\s\p(2m|2n),\o\s\p(r|2s)\times$\\$\o\s\p(2m-r,2(n-s)))$} & $-\frac{1}{2}$ & \makecell{$(r-m)+$\\ $(n-2s)$} & 0 & \makecell{$-2(n-2s)+$\\ $2(m-r)$} & $-\frac{3}{2}$   \\
	\hline
	\makecell{$(\o\s\p(2m+1|2n),$\\$\o\s\p(r|2s)\times$\\ $\o\s\p(2m+1-r,2(n-s)))$} & $-\frac{1}{2}$ & \makecell{$(r-m)+$\\ $(n-2s)-\frac{1}{2}$} & 0 & \makecell{$1-2(n-2s)+$\\ $2(m-r)$} & $-\frac{3}{2}$  \\
	\hline
	$(\o\s\p(2m|2n),\g\l(m|n))$ & $-2$ & 0 & $-\frac{1}{2}$ & 0 & $-\frac{1}{2}$\\
	\hline
	\makecell{$(D(2,1;a),$\\$\o\s\p(2|2)\times\s\o(2))$} & $\alpha$ & 0 & $-\frac{1}{2}$ & 0 & $-\frac{1}{2}$ \\
	\hline
	\makecell{$(\o\s\p(4|2n),$\\$\o\s\p(2|2n)\times\s\o(2)$} & $1$ & 0 & $-\frac{1}{2n}$ & 0 &$-\frac{1}{2n}$\\
	\hline
	$(\a\b(1|3),\s\l(1|4))$ & $-3$ & 0 & $-\frac{5}{4}$ & 0 & $-\frac{1}{4}$ \\
	\hline
	$(\a\b(1|3),\g\o\s\p(2|4))$ & $-\frac{3}{2}$ & 0 & $-\frac{1}{2}$ & 0 & $-\frac{1}{2}$ \\
	\hline
\end{tabular}
\end{center}
\

\

Note that for the third supersymmetric pair we assume $(m,r,s)\neq(2,2,0)$ since this case is special and dealt with later in the table.

As a matter of explanation, the meaning of the parameters is as follows.  In the root system $BC(m,n)$, each real component has three Weyl group orbits determined by the length of the root.  In the first component, the multiplicity $m(\alpha)$ of a short root $\alpha$ is $p$, of the next longest root is $k$, and of the longest root is $q$.  In the second real component, the multiplicity of the short root is $r$, the next longest root $k^{-1}$ and the longest root $s$.  As already stated isotropic roots are required to have multiplicity one.

The deformed bilinear form is given by $B_1+kB_2$, where $B_1,B_2$ are the standard Euclidean inner products on the root system $BC$.  Now each of our restricted root systems can be viewed as $BC(m,n)$ with some multiplicities being set to zero.

\section{Appendix: supersymmetric pairs for $\a\g(1|2)$ and $\a\b(1|3)$} 

In this appendix we give the classification of supersymmetric pairs of the Lie superalgebras $\a\g(1|2)$ and $\a\b(1|3)$, as communicated by V. Serganova.  We refer to Chpt. 26 of \cite{timashev2011homogeneous} for the classification of symmetric pairs of simple Lie algebras.

For $\g=\a\g(1|2),\a\b(1|3)$ all automorphisms are inner by \cite{serganova1985automorphisms}.  Thus we have $\Aut(\a\g(1|2))=SL_2\times G_2$ and $\Aut(\a\b(1|3))=(SL_2\times Spin_7)/\{\pm1\}$. 

In both cases, $\g_{\ol{0}}=\s\l_2\times\k$ for $\k=\g_2$ or $\s\o(7)$.  If $\theta$ is an involution of $\g$ then it is given by $\Ad(g_1g_2)$ where $g_1\in SL_2$ and $g_2\in G_2,Spin_7$, respectively.  Then for $\theta$ to be an involution we must have that $\Ad(g_1^2g_2^2)=\id$, $g_1^2$ is central in $SL_2$, and $g_2^2$ is central in $G_2,Spin_7$ respectively.

The possible choices for $g_1$ up to conjugation are $\pm1$ or $\text{diag}(i,-i)$.  Notice that $\pm1$ induces a trivial involution on $G_0$, while $\text{diag}(i,-i)$ induces a non-trivial involution, but the square of this element is $-1$.  We now do a case by case analysis for each $\g$ for what $g_2$ can be such that we obtain an involution on all of $\g$.  Notice that for $\a\b(1|3)$ we quotient out by $\pm1$, so work up to sign for the choice of $g_2$ in this case.

\textbf{The case $\a\g(1|2)$:} Let us begin with $\a\g(1|2)$.  The center of $G_2$ is trivial, so $g_2$ can either be $1$ or any order two element of $G_2$, of which there is only one up to conjugacy. Thus the only possibility is to have $g_1=\pm1$, i.e.\ we only obtain the two involutions $\Ad(\pm g_2)$, for $g_2$ a fixed non-central element of order 2 in $G_2$.  Further observe that $\Ad(g_2)=\delta\circ\Ad(-g_2)$, and in particular these involutions agree on $\g_{\ol{0}}$.

In each case the Cartan subspace $\a$ is given by a Cartan subalgebra of $\g_2$, and we compute that $\c(\a)=\h+\o\s\p(1|2)$.  Since this has a nontrivial odd part, these involutions $\Ad(g_2)$, $\Ad(-g_2)$ are non-conjugate, and only one satisfies an Iwasawa decomposition.  

Present $\h^*$ and the root systems of $\a\g(1|2)$ as in Section 5.  Then we take $g_2$ to be the element of the maximal torus of $G_2$ which acts by $1$ on $\g_{\epsilon_1}$ and by $-1$ on $\g_{\epsilon_2}$ and $\g_{\epsilon_3}$.  

With this choice of $g_2$, the even roots $\alpha$ for which $\g_{\alpha}$ is fixed by $\Ad(\pm g_2)$ are $\pm 2\delta$, $\pm\epsilon_1$, $\pm(\epsilon_2-\epsilon_3)$.  The odd root spaces fixed by $\Ad(g_2)$ are those with the roots $\pm\delta$, $\pm\delta\pm\epsilon_1$.  Thus the fixed subalgebra of $\Ad(g_2)$ is $\o\s\p(3|2)\times\s\l_2$.  

The odd root spaces fixed by $\Ad(-g_2)$ are $\pm\delta\pm\epsilon_2$, $\pm\delta\pm\epsilon_3$.  Thus the fixed subalgebra of $\Ad(-g_2)$ is $D(2,1;a)$ for some $a$.  To figure out which $a$ we look at the bilinear form on $\h^*$.  We have $(2\delta,2\delta)=-8$, $(\epsilon_1,\epsilon_1)=2$, and $(\epsilon_2-\epsilon_3,\epsilon_2-\epsilon_3)=6$.  Dividing by $2$ we have $a=3$ or $-4$.  This gives that the fixed subalgebra is $D(2,1;3)\cong D(2,1;-4)$.

\textbf{The case $\a\b(1|3)$:} Now we consider $\g=\a\b(1|3)$, and present its root system as in Section 5.  Take $t_1,t_2,t_3$ in the maximal torus of $Spin_7$ to be such that $t_i$ acts on the root space $\g_{\frac{1}{2}(\pm\delta\pm\epsilon_1\pm\epsilon_2\pm\epsilon_3)}$ by $\pm \sqrt{-1}$, where the sign is determined by the sign of $\frac{1}{2}\epsilon_i$ in the root. For example $t_2$ acts by $\sqrt{-1}$ on $\g_{\frac{1}{2}(\pm\delta\pm\epsilon_1+\epsilon_2\pm\epsilon_3)}$ and by $-\sqrt{-1}$ on $\g_{\frac{1}{2}(\pm\delta\pm\epsilon_1-\epsilon_2\pm\epsilon_3)}$. 

Then $t_1,t_2,$ and $t_3$ commute, and we have $t_i^2=-1$ for all $i$. Thus $t_1,t_1t_2,$ and $t_1t_2t_3$ are all square central, and up to $\pm1$ (which we ignore, see the comment towards the beginning of this section), all square central elements of $Spin_7$ are conjugate to one of these three elements, so these are all possibilities we need to consider for $g_2$.  

We observe that $t_1^2=(t_1t_2t_3)^2$ induce multiplication by $(-1)$ on $\g_{\ol{1}}$, while $(t_1t_2)^2=\id$ is the identity on $\g_{\ol{1}}$.   Therefore we obtain involutions of $\g$ given by $\Ad(gt_1)$, $\Ad(t_1t_2)$, and $\Ad(g t_1t_2t_3)$, where $g=\text{diag}(i,-i)\in SL_2$.  These are all of them up to composition with $\delta$.

However we claim that each of these pairs is conjugate to their composition with $\delta$.  Indeed, let $\sigma\in Spin_7$ be the a lift of the element of the Weyl group which sends $\epsilon_1\mapsto-\epsilon_1$ while fixing $\epsilon_2$ and $\epsilon_3.$  Then conjugating one of these involutions by $\sigma$ will have the effect of composing with $\delta$.  In particular, all involutions of $\a\b(1|3)$ admit an Iwasawa decomposition by Corollary \ref{main_cor}.

Now we go through each involution and compute its fixed points.

\textbf{Involution $\Ad(gt_1)$:}  First we look at $\Ad(gt_1)$ which has even fixed root spaces with roots $\pm\epsilon_2\pm\epsilon_3$, $\pm\epsilon_2$, $\pm\epsilon_3$, and odd fixed root spaces with roots $\pm\frac{1}{2}(\delta+\epsilon_1\pm\epsilon_2\pm\epsilon_3)$.  The even part is $\s\o(5)\times\C^2$ and the odd part is the standard representation of $\s\p(4)$ tensored with a sum of two characters for the torus, thus the fixed points subalgebra is $\g\o\s\p(2|4)$.

\textbf{Involution $\Ad(t_1t_2)$:} Next $\Ad(t_1t_2)$ has fixed even root spaces with roots $\pm2\delta$, $\pm\epsilon_1\pm\epsilon_2$, and $\pm\epsilon_3$, and odd fixed root spaces with roots $\pm\frac{1}{2}(\delta\pm(\epsilon_1-\epsilon_2)\pm\epsilon_3)$.  The even fixed subalgebra is $\s\l(2)\times \s\o(4)\times\s\o(3)\cong\s\l_2\times\s\l_2\times\s\l_2\times\s\l_2$, and the odd part is the outer tensor product of the standard representations of three of the copies of $\s\l_2$; the copy of $\s\l_2$ corresponding to $\pm(\epsilon_1+\epsilon_2)$ acts trivially on the odd part.  Thus we find the fixed subalgebra is $\s\l_2\times D(2,1;a)$ for some $a$.  We compute the value of $a$ by looking at the bilinear form.  We see that $(\delta,\delta)=-1$, $(\epsilon_1-\epsilon_2,\epsilon_1-\epsilon_2)=2/3$, and $(\epsilon_3,\epsilon_3)=1/3$.  Thus $a=2$ or $a=-3$. 

\textbf{Involution $\Ad(gt_1t_2t_3)$:}  Finally $\Ad(gt_1t_2t_3)$ has fixed even root spaces $\pm\epsilon_i\pm\epsilon_j$ for $i\neq j$ and fixed odd root spaces $\frac{1}{2}(c_1\delta+c_2\epsilon_1+c_3\epsilon_2+c_4\epsilon_3)$ such that $c_i\in\{\pm1\}$ and $\sum c_i=0$ mod $4$.  Here the even part is $\s\o(6)\times \C$ and the odd part is the spinor rep tensor a character plus the dual spinor rep tensor the dual character.  Thus this is $\s\l(1|4)$.  	


\end{document}